\documentclass[12pt]{article}
\usepackage{epic,cite,amsmath,amsthm,amssymb,texdraw}

\renewcommand{\theequation}{\thesection.\arabic{equation}} 

\newtheorem{thm}{Theorem}[section]
\newtheorem{prop}[thm]{Proposition}

\newcommand{\nc}{\newcommand}

\font\germ=eufm10
\nc{\goth}[1]{{\germ #1}}
\nc{\qbinom}[2]{{\genfrac{[}{]}{0pt}{}{{#1}}{{#2}}}_{q}}
\nc{\bra}[1]{\langle #1 |}
\nc{\ket}[1]{| #1 \rangle}
\nc{\qpf}[1]{({#1}\, ; \, q^4)_{\infty}}
\nc{\br}[1]{\begin{array}{#1}}
\nc{\Ps}[2]{\Phi^{(#1,#2)}}
\nc{\Psz}[2]{\Phi^{(#1,#2)}(\zeta)}
\nc{\Vz}[1]{V^{(#1)}_{\zeta}}
\nc{\wo}[6]{W^I_{#6}\left(\left.
            \begin{array}{ll}{#1}&{#2}\\{#3}&{#4}\end{array}
            \right|{#5}\right)}
\nc{\wt}[6]{W^{II}_{#6}\left(\left.
            \begin{array}{ll}{#1}&{#2}\\{#3}&{#4}\end{array}
            \right|{#5}\right)}
\nc{\wth}[6]{W^*_{#6}\left(\left.
            \begin{array}{ll}{#1}&{#2}\\{#3}&{#4}\end{array}
            \right|{#5}\right)}
\nc{\wec}[5]{\overline{W}^1_{\ell}\left(\left.
            \begin{array}{ll}{#1}&{#2}\\{#3}&{#4}\end{array}
            \right|{#5}\right)}
\nc{\isomo}{\buildrel {\sim} \over \longrightarrow}
\nc{\Czp}{|_0^+}
\nc{\Czm}{|_0^-}
\nc{\Cop}{|_1^+}
\nc{\Com}{|_1^-}
\nc{\Bzp}{\lfloor_0^+}
\nc{\Bzm}{\lfloor_0^-}
\nc{\Bop}{\lfloor_1^+}
\nc{\Bom}{\lfloor_1^-}
\nc{\Tzp}{\lceil_0^+}
\nc{\Tzm}{\lceil_0^-}
\nc{\Top}{\lceil_1^+}
\nc{\Tom}{\lceil_1^-}
\nc{\bpp}{\bullet^{++}}
\nc{\bpm}{\bullet^{+-}}
\nc{\bmp}{\bullet^{-+}}
\nc{\bmm}{\bullet^{--}}
\nc{\bpmpm}{\bullet^{\pm\pm}}
\nc{\bpmmp}{\bullet^{\pm\mp}}
\nc{\Tb}{\lceil}
\nc{\Bb}{\lfloor}
\nc{\ttn}{\kappa}
\nc{\ttm}{\kappa'}
\nc{\tts}{k}
\nc{\ssa}{s}
\nc{\ssb}{t}
\nc{\ssj}{\varepsilon}
\nc{\Aff}{\operatorname{Aff}}
\nc{\Affz}{\operatorname{Aff}(B^{(0)})}
\nc{\Affo}{\operatorname{Aff}(B^{(1)})}
\nc{\calA}{\mathcal{A}}
\nc{\calP}{\mathcal{P}}
\nc{\calR}{\mathcal{R}}
\nc{\calN}{\mathcal{S}}
\nc{\calNo}{\mathcal{N}}
\nc{\calD}{\mathcal{D}}
\nc{\ot}{\otimes}
\nc{\ve}{\varepsilon}
\nc{\be}{\begin{eqnarray}}
\nc{\ee}{\end{eqnarray}}
\nc{\bea}{\begin{eqnarray}}  
\nc{\ena}{\end{eqnarray}}
\nc{\C}{\mathbf{C}}
\nc{\Q}{\mathbf{Q}}
\nc{\R}{\textit{R}}
\nc{\half}{\ensuremath{\frac{1}{2}}}
\nc{\z}{\zeta}
\nc{\La}{\Lambda}
\nc{\s}{\sigma}
\nc{\Hom}{\operatorname{Hom}}
\nc{\End}{\operatorname{End}}
\nc{\vac}{|\textrm{vac}\rangle}
\nc{\dvac}{\langle\textrm{vac}|}
\nc{\id}{\operatorname{id}}
\nc{\la}{\lambda}  
\nc{\bqa}{\begin{eqnarray}}  
\nc{\eqa}{\end{eqnarray}}  
\nc{\ra}{\rightarrow}  
\nc{\lra}{\longrightarrow}
\nc{\uqp}{U^{\prime}_q (\widehat{sl}_2)}
\nc{\uq}{U_q(\slth)}
\nc{\vsl}{V(\sigma(\lambda))}
\nc{\Rbar}{\bar{R}}
\nc{\vl}{V(\lambda)}  
\nc{\bu}{\bullet}
\nc{\an}{{\ell}}
\nc{\ak}{{\ell}}
\nc{\am}{{\ell'}}
\nc{\al}{{\ell'}}
\nc{\ep}{\varepsilon} 
\nc{\de}{\delta}  
\nc{\slth}{\widehat{\goth{sl}}_2\hskip 1pt}
\nc{\ws}{\;\;}
\nc{\qu}{{1\ov 4}}
\nc{\tf}{\tilde{f}}
\nc{\te}{\tilde{e}}
\nc{\ez}{|_{0}}
\nc{\eo}{|_{1}}
\nc{\cz}{[_{0}}
\nc{\co}{[_{1}}
\nc{\hif}{\hb{ if }}
\nc{\hev}{\hb{ is even }}
\nc{\hod}{\hb{ is odd }}
\nc{\HR}{{_h}\bar{R}}
\nc{\PR}{\bar{R}}
\nc{\cP}{\mathcal{P}}
\nc{\er}{\end{array}}
\nc{\Tr}{{\rm Tr}}
\nc{\hb}{\hbox}
\nc{\cT}{\mathcal{T}}
\nc{\zi}{\zeta^{-1}}
\nc{\nn}{\nonumber} 
\nc{\Z}{\mathbf{Z}}
\nc{\bR}{\bar{R}}
\nc{\curlra}{\buildrel{\sim}\over\longrightarrow}
\nc{\bs}{\tilde{s}}
\nc{\cD}{\mathcal{D}}  
\nc{\cH}{\mathcal{H}}
\nc{\cF}{\mathcal{F}}
\nc{\cO}{\mathcal{O}}
\nc{\epp}{\varepsilon^{\prime}} 
\nc{\bi}{\bar{i}}
\nc{\bj}{\bar{j}}
\nc{\ol}{\overline}
\nc{\Om}{\Omega}
\nc{\pl}{\prod\limits} 
\nc{\sli}{\sum\limits} 

\begin{document}

\begin{center}
\vspace*{10mm}
{\Large \bf Determinant Formula for the Solutions
of the Quantum Knizhnik-Zamolodchikov Equation with $|q|=1$\\[8mm]}
{Tetsuji Miwa and Yoshihiro Takeyama\\[8mm]
Research Institute for Mathematical Sciences,
Kyoto University, Kyoto 606, Japan\\[10mm]
{\it Dedicated to Professor Kazuhiko Aomoto on his sixtieth birthday.}
}
\end{center}

\begin{abstract}
\noindent 
The fundamental matrix solution of the quantum Knizhnik-Zamolodchikov
equation associated with $U_q(\widehat{sl}_2)$ is constructed for $|q|=1$.
The formula for its determinant is given in terms of the double
sine function.
\end{abstract}
\nopagebreak

\setcounter{equation}{0}
\makeatletter
\def\eqnarray{\stepcounter{equation}\let\@currentlabel\theequation
\global\@eqnswtrue\m@th
\global\@eqcnt\z@\tabskip\@centering\let\\\@eqncr
$$\halign to\displaywidth\bgroup\@eqnsel\hskip\@centering
  $\displaystyle\tabskip\z@{##}$&\global\@eqcnt\@ne
   \hfil${{}##{}}$\hfil
  &\global\@eqcnt\tw@ $\displaystyle\tabskip\z@{##}$\hfil
   \tabskip\@centering&\llap{##}\tabskip\z@\cr}
\makeatother
\def\bea{\begin{eqnarray}}
\def\ena{\end{eqnarray}}
\def\no{\nonumber}
\def\sh{\mathop{\rm sh}\nolimits}
\def\ch{\mathop{\rm ch}\nolimits}
\def\Re{\mathop{\rm Re}\nolimits}
\def\Im{\mathop{\rm Im}\nolimits}
\def\qKZ{\hbox{$q${\rm KZ }}}
\section{Introduction} 
The aim of this paper is to compute the determinant of
the fundamental matrix solution of the quantum Knizhnik-Zamolodchikov equation
associated with $U_q(\widehat{sl}_2)$ for $|q|=1$.

In \cite{TV1} and
\cite{TV2}, Tarasov and Varchenko studied the rational \qKZ equation
associated with the $sl_2$ Yangian and the trigonometric \qKZ equation
associated with $U_q(\widehat{sl}_2)$, respectively. The case studied in
\cite{TV1} is the limit $q\to1$ of the case studied in \cite{TV2}.
However, there is a difference in these two cases. In the former,
the unknown function $\psi(\beta_1,\ldots,\beta_n)$
is defined on ${\bf C}^n$, while in the latter
the unknown function $\tilde\psi(z_1,\ldots,z_n)$
is defined on $({\bf C}\backslash\{0\})^n$. If we set
\bea
q=e^{-\frac{\pi^2i}\rho},\quad z_m&=&e^{\frac{2\pi\beta_m}\rho}(1\leq m\leq n),
\no
\ena
the limit $q\to1$ corresponds to $\rho\to\infty$. Therefore,
the period $\rho i$ in $\beta_m$ is lost in this limit.

We will study the \qKZ equation on ${\bf C}^n$ with finite $\rho$.
In \cite{TV2}, the solutions are single-valued on
$({\bf C}\backslash\{0\})^n$. The price for this restriction is that the
multiplicative step $p$ in $z_m$ of the difference equation is restricted
by the condition $0<|p|<1$ (or $|p|>1$ in a different convention).
We consider the case
$p=e^{-\frac{2\pi i\lambda}\rho}$ with real $\rho$ and $\lambda$.
Namely, the additive step in $\beta_m$ is $-\lambda i$.
Then, the solutions are {\it not} single-valued in $z_m$.

In the applications to physics, this difference corresponds to the
difference of the models. In the application to the XXZ spin
chain in the massive regime \cite{JM}, the deformation parameter $q$ satisfies
$-1<q<0$, and the multiplicative step is given by $p=q^{-4}$ for the
correlation functions and by $p=q^4$ for the form factors.
On the other hand, in the application to the sine-Gordon model \cite{Smirnov},
the deformation parameter $q$ satisfies $|q|=1$ and the additive step is given
by $\lambda=2\pi$ (i.e., $p=q^4$).

We follow \cite{M1,V,TV2} for the construction of solutions of the \qKZ
equations. We consider the evaluation modules
of the $U_q(sl_2)$ Verma modules with the spectral parameter $z_m$
and the highest weight $2\Lambda_m$. The unknown function
$\psi(\beta_1,\ldots,\beta_n)$ takes a value in their tensor product.
We will fix $l$ and consider the subspace of weight
$2(\sum_{m=1}^n\Lambda_m-l)$, which has dimensions $d=\binom{n+l-1}{n-1}$.
The \qKZ equation is defined on this subspace. We will choose a basis of
this space. The \qKZ equation is then written as a system of equations
for a $d$-dimensional unknown function.

We will construct a non-degenerate $d\times d$ matrix such that each column
solves the \qKZ equation following a similar construction in \cite{TV1}.
We call it the fundamental solution.
The fundamental solution is constructed symmetrically with respect to the
parameters $\rho$ and $\lambda$. Namely, each row solves another \qKZ
equation with
\bea
q=e^{-\frac{\pi^2i}\lambda},\quad z_m&=&e^{\frac{2\pi\beta_m}\lambda}
(1\leq m\leq n),
\no
\ena
and $-\rho i$ being the additive step in $\beta_m$.
We remark that in \cite{FTV} the co-existence of two systems of equations
has already been discussed in the setting of elliptic equations.

We give the fundamental solution in terms of $l$-dimensional integrals.
To prove the non-degeneracy of the matrix we compute its determinant.
The integrands of these integrals are meromorphic functions given
explicitly in terms of the double sine function.
We express the determinant itself as a simple product
of double sine functions. In \cite{TV1,TV2}, similar determinant formulas
are obtained. We follow their line of arguments in the evaluation
of the determinant.

We introduce another parameter $\mu$
which couples to the $sl_2$ generator $h_1$ acting on the $m$-th component
of the tensor product. Therefore, the determinant is a function of the
variables $\beta_1,\ldots,\beta_n,\Lambda_1,\ldots,\Lambda_n$
and $\rho,\lambda,\mu$. The dependence on $\beta_1,\ldots,\beta_n$
is easily determined by using the two sets of \qKZ equations.
In order to determine the multiplicative factor which is independent of
$\beta_1,\ldots,\beta_n$, we compute the asymptotics of the fundamental
solution in the region $\beta_1<\hskip-2pt<\cdots<\hskip-2pt<\beta_n$.
This asymptotics is triangular, and the diagonal extries are the same kind
of integrals with $n=1$.

The computation of the $n=1$ integrals is done in two steps. The first step
is to derive difference equations in $\mu$ with steps $\frac{2\pi}\rho$
and $\frac{2\pi}\lambda$. The dependence on $\mu$ is determined by these 
equations. Finally, the multiplicative factor which is independent of
$\mu$ is determined by taking the asymptotics when $\mu\to i\infty$.

The plan of the paper is as follows. Section 2 contains preliminaries
of the \qKZ equation. In Section 3, we construct the fundamental matrix
solution. In Section 4, the asymptotics of
the determinant is computed in terms of the $n=1$ integrals. These integrals
are computed in Section 5. The formula for the determinant is given 
in Section 6. A short Appendix is given on the double sine function for the
convenience of the reader.

\section{The \qKZ equation}

Let $q$ be a nonzero complex number which is not a root of unity.
Consider the quantum affine algebra $\widehat{U_{q}}=U'_{q}(\widehat{sl_{2}})$
with generators $e_i, f_i, q^{\pm h_i} \, (i=0,1)$ and relations:
\bea 
q^{h_i}q^{h_j}=q^{h_j}q^{h_i}, & & q^{h_i}q^{-h_i}=q^{-h_i}q^{h_i}=1, \no \\
q^{h_i}e_{j}q^{-h_i}=q^{a_{ij}}e_{j}, 
& & q^{h_i}f_{j}q^{-h_i}=q^{-a_{ij}}f_{j}, \no \\
\left[ e_{i} , f_{j} \right] &=& \displaystyle{\delta_{ij} \frac{q^{h_i}-q^{-h_i}}{q-q^{-1}}},  \no\\
e_i^{3}e_{j}-[3]_{q}e_i^{2}e_{j}e_{i} &+& [3]_{q}e_{i}e_{j}e_{i}^{2}-e_{j}e_{i}^{3}=0 \quad (i \not= j), \no \\
f_i^{3}f_{j}-[3]_{q}f_i^{2}f_{j}f_{i} &+& [3]_{q}f_{i}f_{j}f_{i}^{2}-f_{j}f_{i}^{3}=0 \quad (i \not= j). \no
\ena
Here $a_{00}=a_{11}=2, a_{01}=a_{10}=-2,$ and $[n]_{q}$ is the $q$-integer
defined by $[n]_{q}=\frac{q^{n}-q^{-n}}{q-q^{-1}}$. 
We also use $[n]_{q}!=\prod_{m=1}^n[m]_{q}$.

We use the coproduct
\bea
\Delta (q^{h_i})=q^{h_i} \otimes q^{h_i}, \quad \Delta (e_i)=e_i \otimes 1 +
q^{h_i} \otimes e_i, \quad \Delta (f_i)=f_i \otimes q^{-h_i} + 1 \otimes f_i.
\no
\ena

We denote by $U_q$ the subalgebra generated by $e_1, f_1, q^{\pm h_1}.$
Let $\Lambda$ be a complex number and $V_{\Lambda}^q$ be the Verma module
for $U_q$ with the highest weight $q^{2\Lambda}$.
We denote the highest weight vector of $V^q_\Lambda$ by $v^{(0)}$
and use the following basis $\{v^{(k)}\}_{k\geq0}$.
\bea 
\quad e_1 v^{(k)}=[2\Lambda-k+1]_qv^{(k-1)},
\quad q^{h_1}v^{(k)}=q^{2(\Lambda-k)}v^{(k)},
\quad f_1v^{(k)}=[k+1]_qv^{(k+1)}.\no
\ena
For a nonzero complex number $z$, we can give $V_{\Lambda}^q $ the $\widehat{U_q}$ module structure by setting 
\bea 
e_0v=zf_1v,\quad q^{h_0}v=q^{-h_1}v,\quad f_0v=z^{-1}e_1v,\quad
(v\in V^q_\Lambda).
\no
\ena
We denote by $V_{\Lambda}^q (z)$ the $\widehat{U_q}$ module defined in this
way.

Let $\Lambda_1, \Lambda_2$ be complex numbers, and $v_1^{(k)},v_2^{(k)}$
the base vectors of $V_{\Lambda_1}^q,V_{\Lambda_2}^q$,
respectively. For generic $z_1,z_2$, there exists a map
$\displaystyle R_{\Lambda_{1}\Lambda_{2}}^q(z_1/z_2) \in {\rm End}
\left(V_{\Lambda_1}^q(z_1)\otimes V_{\Lambda_2}^q(z_2)\right)$ such that
\bea 
R_{\Lambda_{1}\Lambda_{2}}^{\, q}(z_1/z_2)v_1^{(0)} &\otimes&  v_2^{(0)} = v_1^{(0)} \otimes v_2^{(0)}, \label{eq:norm} \\
R_{\Lambda_{1}\Lambda_{2}}^{\, q}(z_1/z_2) \circ \Delta (x) &=& \sigma ( \Delta (x)) \circ R_{\Lambda_{1}\Lambda_{2}}^{\, q}(z_1/z_2) \quad {\rm for \, all} \, x \in \widehat{U_q},
\label{eq:compa} \ena
where $\sigma$ is the permutation map $\sigma ( a \otimes b )= b \otimes a$. Such a map  $\displaystyle R_{\Lambda_{1}\Lambda_{2}}^{\, q}(z_1/z_2) $ is uniquely determined by (\ref{eq:norm}) and (\ref{eq:compa}).

Let $V_{\Lambda_1}^q\otimes V_{\Lambda_2}^q=\oplus_{l=0}^{\infty}V_{\Lambda_1
+\Lambda_2-2l}^{\,q}$
be the decomposition of the $U_q$ module into the irreducible components. 
We denote by $\Pi_l$ the projection to the component
$V_{\Lambda_1+\Lambda_2-2l}^{\, q}$.
Then we have (see, e.g., \cite{TV1})
\bea
R_{\Lambda_{1}\Lambda_{2}}^{\,q}(z_1/z_2)=R_{\Lambda_{1}\Lambda_{2}}^{\,q}(0)
\displaystyle\sum_{l=0}^{\infty}\Pi_{l}\prod_{j=0}^{l-1}
\frac{1-q^{-2(\Lambda_1+\Lambda_2-j)}z_1/z_2}
{1-q^{2(\Lambda_1+\Lambda_2-j)}z_1/z_2},\label{eq:Rmat}
\ena
where
\bea
R_{\Lambda_{1}\Lambda_{2}}^{\,q}(0)=q^{2\Lambda_1\Lambda_2
-\frac{h_1\otimes h_1}{2}}\displaystyle\sum_{k=0}^{\infty}
\frac{q^{k}(1-q^2)^{2k}}{\prod_{j=1}^{k}(1-q^{2j})}(q^{-h_1}
e_1\otimes q^{h_1}f_1)^k.
\no
\ena
Note that $R_{\Lambda_{1}\Lambda_{2}}^{\,q}(0)$ is triangular
with respect to the basis $v^{(j)}_1\otimes v^{(k)}_2$.

Consider the tensor product
$W^q =V_{\Lambda_1}^q (z_1)\otimes \cdots \otimes V_{\Lambda_n}^q (z_n)$.
Let $\displaystyle R_{ij}^{\, q}(z_i/z_j) \in {\rm End}W^q $ be the operator
acting on the $i$-th and $j$-th components as 
$R_{\Lambda_{i},\Lambda_{j}}^{q}(z_i/z_j)\in{\rm End}(V_{\Lambda_{i}}^{q}(z_i)
\otimes V_{\Lambda_{j}}^{q}(z_j))$.
We denote by $H_m$ the operator on $W^q $ given by
\bea
H_m=1\otimes\cdots\otimes\stackrel{m-{\rm th}}{h_1}\otimes\cdots\otimes1.
\no
\ena

Let $p, r$ be non-zero complex numbers. For $m=1, \ldots , n$, set
\bea
\tilde K_{m}(z_1, \cdots , z_n;q,p,r) &=& R_{m,m-1}^{\, q}(pz_m/z_{m-1})
\cdots R_{m,1}^{\, q}(pz_m/z_{1}) r^{2\Lambda_{m}-H_m}\no\\
&&{}\times R_{m,n}^{\,q}(z_m/z_n)\cdots R_{m,m+1}^{\,q}(z_m/z_{m+1}).
\no
\ena
The quantum Knizhnik-Zamolodchikov equation is the following system
of equations for unknown function $\tilde\psi(z_1, \cdots ,z_n)$
that takes a value in $W^q $:
\bea
\tilde \psi(z_1, \cdots , pz_m, \cdots , z_n)=
\tilde K_{m}(z_1, \cdots , z_n ; q, p, r)
\tilde \psi(z_1, \cdots , z_n),
(1\leq m\leq n)\label{eq:qkz}
\ena

In this paper, we consider the \qKZ equation where the step $p$ satisfies
$|p|=1$.
In this case, solutions to (\ref{eq:qkz}) are not necessarily single-valued
with respect to the parameters $z_1,\ldots,z_n,q,p,r$. Therefore,
we rewrite (\ref{eq:qkz}) as follows. We set
\bea
&&z^{(\rho)}_k=e^{\frac{2\pi\beta_k}{\rho}},q^{(\rho)}=e^{-\frac{\pi^2 i}{\rho}},
p^{(\rho,\lambda)}=e^{-\frac{2\pi i\lambda}{\rho}},
r^{(\lambda)}=e^{-\frac{\mu\lambda i}{2}},
\no\\
&&\psi(\beta_1,\cdots,\beta_n)=\tilde\psi(z^{(\rho)}_1,\cdots,z^{(\rho)}_n),\no\\
&&K_{m}(\beta_1,\cdots,\beta_n;\rho,\lambda,\mu)=
\tilde K_{m}(z^{(\rho)}_1,\cdots,z^{(\rho)}_n;q^{(\rho)},p^{(\rho,\lambda)},
r^{(\lambda)}).
\label{eq:rewrite}
\ena
Then (\ref{eq:qkz}) is equivalent to the following system of equations
for $\psi$:
\bea
\psi(\beta_1 , \cdots , \beta_m -\lambda i , \cdots ,\beta_n)
=K_{m}(\beta_1,\cdots,\beta_n;\rho,\lambda,\mu)
\psi(\beta_1 , \cdots , \beta_n),
(1\leq m\leq n).\label{eq:qkz1}
\ena

For a non-negative integer $l$, set
\bea
W^q_l=\{v\in W^q;h_1v=2(\sum_{m=1}^n\Lambda_m-l)v\}.
\no
\ena
Since the matrix $K_m$ acts on $W^q_l$ for each $l$, the equation (\ref{eq:qkz1})
for $\psi$ splits into the equations for the weight-$l$ component
$\psi_l\in W^q_l$.
We denote by $K_{m,l}$ the $W^q_l$ block of $K_m$.
Then, we have
\bea
\psi_l(\beta_1,\cdots,\beta_m-\lambda i , \cdots ,\beta_n)
=K_{m,l}(\beta_1,\cdots,\beta_n;\rho,\lambda,\mu)
\psi_l(\beta_1 , \cdots , \beta_n),
(1\leq m\leq n).\label{eq:qkz2}
\ena

Let $W^{q*}_l$ be the dual space of $W^q_l$ with the coupling
$\langle v^*,v\rangle$ for $v^*\in W^{q*}_l$ and $v\in W^q_l$.
In Section 3, we construct $\Psi_l(\beta,\ldots,\beta_n;\rho,\lambda,\mu)$
taking a value in $W^{q^{(\rho)}}_l\otimes W^{q^{(\lambda)}}_l$
such that for any $v^*\in W^{q^{(\lambda)}*}_l$,
\bea
\psi_l(\beta_1,\ldots,\beta_n)&=&
\langle v^*,\Psi_l(\beta,\ldots,\beta_n;\rho,\lambda,\mu)\rangle
\in W^{q^{(\rho)}}_l\no
\ena
satisfies (\ref{eq:qkz2}) and also, for any $w^*\in W^{q^{(\rho)*}}_l$,
\bea
\psi_l(\beta_1,\ldots,\beta_n)&=&
\langle w^*,\Psi_l(\beta,\ldots,\beta_n;\rho,\lambda,\mu)\rangle
\in W^{q^{(\lambda)}}_l\no
\ena
satisfies (\ref{eq:qkz2}) with $\rho$ and $\lambda$ interchanged.

\section{The hypergeometric solutions}

For a non-negative integer $l$, we set
\bea 
{\cal Z}_{l}^{n}=\left\{ (l_1, \cdots , l_n)
\in {\bf Z}_{\ge 0}^{n};\sum_{m=1}^{n} l_m = l \right\}.
\no
\ena
For $L=(l_1, \cdots , l_n) \in {\cal Z}_{l}^{n}$,
set
\bea
\Gamma_{m}^{L}=\{l_1+\cdots+l_{m-1}+1, \cdots ,
l_1+\cdots+l_m\},
\no
\ena
and define a map $\gamma^{L}:\{1,\cdots,l\}\to\{1,\cdots,n\}$ by
\bea 
\gamma^{L}(j)=m \stackrel{{\rm def}}{\Longleftrightarrow} j \in \Gamma_{m}^{L}.
\no
\ena 
We also define a partial order $\preceq$ on ${\cal Z}_{l}^{n}$ by
\bea
L \preceq L' \stackrel{{\rm def}}{\Longleftrightarrow}
l_n\leq l'_n,l_{n-1}+l_n\leq l'_{n-1}+l'_n,\ldots,l_2+\cdots+l_n\leq
l'_2+\cdots+l'_n.
\no
\ena

We define a function $w_{L}^{(\rho)}$ as follows:
\bea
w_{L}^{(\rho)}(\alpha_1,\cdots,\alpha_{l};\beta_{1},\cdots,\beta_{n})
&=&{\rm Skew}\left(g_{L}^{(\rho)}(\alpha_1,\cdots,\alpha_{l};
\beta_{1},\cdots,\beta_{n})\right),\\
g_{L}^{(\rho)}(\alpha_1,\cdots,\alpha_{l};\beta_{1},\cdots,
\beta_{n})&=&\displaystyle{q^{(\rho)}}^{\sum_{m<m'}l_m l_{m'}} 
\prod_{1 \le j < j' \le l}\sh{\frac{\pi}{\rho}(\alpha_{j'} 
- \alpha_j - \pi i)}\no\\
\times \prod_{j=1}^{l}
\Biggl( e^{-\frac{\pi}{\rho}(\alpha_j -\beta_{\gamma^L(j)}+\Lambda_{\gamma^L(j)} \pi i)}
&\displaystyle \prod_{m<{\gamma^L(j)}}& \left.  \sh{\frac{\pi}{\rho}
(\alpha_j -\beta_m +\Lambda_m \pi i)} \prod_{m>{\gamma^L(j)}} 
\sh{\frac{\pi}{\rho}(\alpha_j -\beta_m -\Lambda_m \pi i)} \right),
\label{eq:def1}
\no
\ena
where Skew is the skew-symmetrization with respect to $(\alpha_1, \cdots , \alpha_{l})$, i.e.,
\bea
{\rm Skew}f(\alpha_1, \cdots , \alpha_{l})=\frac{1}{l !}\sum_{\sigma \in {\bf S}_{l}} ({\rm sgn}\sigma) f(\alpha_{\sigma_1}, \cdots , \alpha_{\sigma_{l}}).
\no
\ena
We abbreviate 
$w_{L}^{(\rho)}(\alpha_1,\cdots,\alpha_{l};\beta_{1},\cdots,\beta_{n})$
to
$w_{L}^{(\rho)}(\beta_{1},\cdots,\beta_{n})$
(or $w_{L}^{(\rho)}(\alpha_1,\cdots,\alpha_{l})$)
when the dependence on the abbreviated variables is irrelevant.

We set
\bea{\cal F}^{(\rho)}=\sum_{L \in {\cal Z}_{l}^{n}}{\bf C}w_{L}^{(\rho)}.
\no
\ena
Now we define a pairing between ${\cal F}^{(\rho)}$ and ${\cal F}^{(\lambda)}$.
We use
\bea
\varphi(x ;\Lambda)=\frac{1}{S_{2}(ix -\Lambda \pi)S_{2}(-ix -\Lambda \pi)},
\quad\phi(x)=\prod_{m=1}^{n}\varphi(x -\beta_{m};\Lambda_{m}),
\quad\psi(x)=\varphi(x;-1),\no\\ \label{eq:defphi}
\ena
where $S_{2}(x)=S_{2}(x|\rho,\lambda)$ is the double sine function
with periods $\rho$ and $\lambda$ (see Appendix).

In this paper, we assume that the auxiliary parameters,
$\rho$, $\lambda$, $\mu$ and $-\Lambda_m$ $(1\leq m\leq n)$, are
in ${\bf R}_{>0}$.

Suppose that $f,g$ are entire functions in the variables
$\alpha_1,\ldots,\alpha_l$. For fixed $\beta_1,\ldots,\beta_n\in{\bf R}$
we set
\bea
I(f,g)&=&\int_{C} \cdots \int_{C}
d\alpha_{1} \cdots d\alpha_{l} \, e^{\mu \sum_{j=1}^{l} \alpha_{j}}
\prod_{j=1}^{l} \phi (\alpha_{j}) \! \! \prod_{1 \le j < j' \le l}
\! \! \psi(\alpha_{j}-\alpha_{j'})\nonumber\\
&\times&f(\alpha_1,\ldots,\alpha_l)g(\alpha_1,\ldots,\alpha_l).\label{eq:hyp}
\ena
Here the contour $C$ for the variable $\alpha_j $ is taken to be
the real line ${\bf R}$.
Note that the poles of the integrand at
\bea
\beta_{m}-\Lambda_{m}\pi i+\rho i{\bf Z}_{\ge0}+\lambda i{\bf Z}_{\ge 0},
\quad\alpha_{j'}+\pi i+\rho i{\bf Z}_{\ge0}+\lambda i{\bf Z}_{\ge0}\no
\ena
are above $C$ and the poles at 
\bea
\beta_{m}+\Lambda_{m}\pi i-\rho i{\bf Z}_{\ge0}-\lambda i{\bf Z}_{\ge0},
\quad\alpha_{j'}-\pi i-\rho i{\bf Z}_{\ge0}-\lambda i{\bf Z}_{\ge0}\no
\ena
are below $C$. As we will see in Section 4, there is a region of the parameters
where the integral (\ref{eq:hyp}) is absolutely convergent.

The integral, in particular, defines the pairing $I(w^{(\rho)},w^{(\lambda)})$
between $w^{(\rho)}\in{\cal F}^{(\rho)}$ and
$w^{(\lambda)}\in {\cal F}^{(\lambda)}$ .

\begin{thm}
For $w^{(\lambda)} \in {\cal F}^{(\lambda)}$, we set
\bea
\psi_l&=&\sum_{L \in {\cal Z}_{l}^{n}}I(w_{L}^{(\rho)},w^{(\lambda)})
v^{(\rho)}_L\hbox{ where }
v^{(\rho)}_L=v_{1}^{(l_1)}\otimes\cdots\otimes v_{n}^{(l_n)}\in V_{\Lambda_1}
^{q^{(\rho)}}(z^{(\rho)}_1)\otimes\cdots\otimes
V_{\Lambda_n}^{q^{(\rho)}}(z^{(\rho)}_n).
\no
\ena
Then, $\psi_l$ is a solution to $(\ref{eq:qkz2})$.
\end{thm}

\begin{proof}
For $L=(l_1,\ldots,l_n)$, set
$L^{(m)}=(l_1,\cdots,l_{m+1},l_{m},\cdots,l_{n})$ and
${\overline L}=(l_n , l_1 , \cdots , l_{n-1})$.
It is easy to see that $\psi_l$ is a solution to
(\ref{eq:qkz2}) if the following relations (\ref{eq:rel1}) and
(\ref{eq:rel2}) hold.
\bea
&&\sum_{L\in{\cal F}^n_l}
w_{L^{(m)}}^{(\rho)}(\beta_{1},\cdots,\beta_{m+1},\beta_{m},\cdots,\beta_{n})
v^{(\rho)}_L\no\\
&&=\sum_{L\in{\cal F}^n_l}
R^{q^{(\rho)}}_{\Lambda_m,\Lambda_{m+1}}(z^{(\rho)}_{m}/z^{(\rho)}_{m+1})
w^{(\rho)}_L(\beta_{1},\cdots,\beta_{n})v^{(\rho)}_L,
\label{eq:rel1}
\ena
\bea
I(w_{{\overline L}}^{(\rho)}
(\beta_n-\lambda i,\beta_1,\cdots,\beta_{n-1}),
w^{(\lambda)}(\beta_1,\cdots,\beta_n-\lambda i))&=&e^{-\mu l_n\lambda i}
I(w_{L}^{(\rho)},w^{(\lambda)}).
\label{eq:rel2}
\ena

The relation (\ref{eq:rel1}) follows from Lemma 5.2.2 in \cite{M1}. We prove
(\ref{eq:rel2}) in the form
\bea
I(g_{L}^{(\rho)} , w^{(\lambda)})(\beta_1 , \cdots , \beta_n+\lambda i)
&=&e^{\mu l_n \lambda i}I(g_{{\overline L}}^{(\rho)}
(\beta_n, \beta_1, \cdots , \beta_{n-1}), w^{(\lambda)}
(\beta_1 , \cdots , \beta_n )).\label{eq:rel3}
\ena

First note that 
\bea
w^{(\lambda)}(\alpha_1,\cdots,\alpha_j+\lambda i,\cdots,\alpha_{l})
&=&(-1)^{n+l-1}w^{(\lambda)}(\alpha_1,\cdots,\alpha_j,\cdots,\alpha_{l}),
\label{eq:causep1}\\
w^{(\lambda)}(\beta_1,\cdots,\beta_{n-1},\beta_{n}+\lambda i)
&=&(-1)^{l}w^{(\lambda)}(\beta_1,\cdots,\beta_{n-1},\beta_{n}).
\label{eq:causep2}
\ena
From (\ref{eq:Sprop}), we find
\bea
\frac{\varphi(x-\lambda i;\Lambda)}{\varphi(x;\Lambda)}
=-\frac{\sh{\frac{\pi}{\rho}(x+\Lambda\pi i)}}
{\sh{\frac{\pi}{\rho}(x-\Lambda\pi i-\lambda i)}}.\label{eq:causep3}
\ena

Set
\bea
f_{L}(\alpha_1,\cdots,\alpha_{l};\beta_1,\cdots,\beta_n)
&=&e^{\mu\sum_{j=1}^{l}\alpha_j}\prod_{j=1}^{l}
\phi(\alpha_j)\prod_{1\le j<j'\le l}
\psi(\alpha_j-\alpha_{j'})\no\\
&&\times g_{L}^{(\rho)}(\alpha_1,\cdots,\alpha_{l};\beta_1,\cdots,\beta_n).
\no
\ena
Then
\bea
&&I(g_{L}^{(\rho)},w^{(\lambda)})(\beta_1,\cdots,\beta_n)\no\\
&&=\int_{C}\cdots\int_{C}d\alpha_1\cdots d\alpha_{l}
f_{L}(\alpha_1,\cdots,\alpha_{l};\beta_1,\cdots,\beta_n)
w^{(\lambda)}(\alpha_1,\cdots,\alpha_{l};\beta_1,\cdots,\beta_n).
\no
\ena

First, we consider the case $l_n =0$. In this case, $f_L$ has no poles
at
\bea
\alpha_j&=&\beta_n+\Lambda_n\pi i-\rho i{\bf Z}_{\ge0}.\no
\ena
Therefore, when we make the
analytic continuation $\beta_n \to \beta_n +\lambda i$, no poles of $f_{L}$
cross the contour $C$. From (\ref{eq:causep2}) and (\ref{eq:causep3}), we find
\bea
&&f_{L}(\beta_1,\cdots,\beta_n+\lambda i)
w^{(\lambda)}(\beta_1,\cdots,\beta_n+\lambda i)
=f_{{\overline L}}(\beta_n,\beta_1,\cdots,\beta_{n-1})
w^{(\lambda)}(\beta_1,\cdots,\beta_n).\no
\ena
Therefore, the equation (\ref{eq:rel3}) holds if $l_n =0$.

In the case of $l_n>0$, if $j\in\Gamma^L_n$, the poles of $f_{L}$ at
\bea
\alpha_j=\beta_n+\Lambda_n\pi i-\rho i{\bf Z}_{\ge0}
\no
\ena
may cross the contour $C$. In order to avoid this crossing,
we shift the contour $C$ to $C+\lambda i$ for all $j\in\Gamma^L_n$.
We note that $f_{L}$ has no poles at
\bea
\alpha_j&=&\beta_m-\Lambda_m\pi i+\rho i{\bf Z}_{\ge0}\quad
(m=1,\cdots,n-1),\no\\
\alpha_j&=&\alpha_{j'}+\pi i+\rho i{\bf Z}_{\ge0}\quad
(j'\in\Gamma_{m}^{L},m=1,\cdots n-1)\no
\ena
for $j \in \Gamma_{n}^{L}$. Hence, this shift of the contours
does not cause any crossing of poles. 

By changing the variables $\alpha_j \to \alpha_j + \lambda i \,
(j \in \Gamma_{n}^{L})$ after this shift and using (\ref{eq:causep1}),
(\ref{eq:causep2}), and (\ref{eq:causep3}), we get
\bea
&&I(g_{L}^{(\rho)},w^{(\lambda)})(\beta_1,\cdots,\beta_n+\lambda i)
=(-1)^{l_n(l-l_n)}e^{\mu l_n\lambda i}\int_{C}\cdots\int_{C}
d\alpha_1\cdots d\alpha_{l}\no\\
&&\quad\times f_{{\overline L}}(\alpha_{l-l_n+1},\cdots,\alpha_{l},
\alpha_1,\cdots,\alpha_{l-l_n};\beta_n,\beta_1,\cdots,
\beta_{n-1})
w^{(\lambda)}(\alpha_1,\cdots,\alpha_{l};\beta_1,\cdots,\beta_{n}).
\label{eq:causep4}
\ena
Symmetrizing the integrand of (\ref{eq:causep4}), we obtain (\ref{eq:rel3}). 
\end{proof}

We define the fundamental matrix solution by
\bea
\Psi_l(\beta_1,\ldots,\beta_n;\rho,\lambda,\mu)
=\sum_{L,L'\in{\cal Z}^n_l}I(w^{(\rho)}_L,w^{(\lambda)}_{L'})
v^{(\rho)}_L\otimes v^{(\lambda)}_{L'}\in W^{(\rho)}_l\otimes W^{(\lambda)}_l.
\ena
Then, it has the property announced in Section 2.

\section{Asymptotics of the Determinant in $\beta_1 , \cdots , \beta_n $} 

In the following sections we calculate the determinant
\bea
D_{l}(\beta_1,\cdots,\beta_n)=\det{\left[I(w_{L}^{(\rho)},w_{L'}^{(\lambda)})
\right]_{L,L'\in{\cal Z}_{l}^{n}}}.
\no
\ena

From Theorem 3.1, we find
\bea
\frac{D_{l}(\beta_1, \cdots , \beta_{m}-\lambda i, \cdots , \beta_{n})}
{D_{l}(\beta_{1}, \cdots ,\beta_{m}, \cdots , \beta_{n})}
&=& {\det}K_{m,l}
(\beta_1, \cdots , \beta_n ;\rho,\lambda,\mu ), \label{eq:1de} \\
\frac{D_{l}(\beta_1, \cdots , \beta_{m}-\rho i, \cdots , \beta_{n})}
{D_{l}(\beta_{1}, \cdots ,\beta_{m}, \cdots , \beta_{n})} 
&=& {\det}K_{m,l}
(\beta_1, \cdots , \beta_n ;\lambda,\rho,\mu ). \label{eq:2de} 
\ena
Using the formula (\ref{eq:Rmat}), we see that
\bea
&&{\det}K_{m,l}(\beta_1, \cdots , \beta_n;\rho,\lambda,\mu)\no\\
&=& \left( e^{- \mu \lambda i} \right)^{\binom{n+l -1}n}
\prod_{j=0}^{l -1} \left( \prod_{k=1}^{m-1} \frac{\sh{\frac{\pi}
{\rho}(\beta_{m}-\beta_{k}-\lambda i+(\Lambda_{m}+\Lambda_{k}-j)\pi i)}}
{\sh{\frac{\pi}{\rho}(\beta_{m}-\beta_{k}-\lambda i
-(\Lambda_{m}+\Lambda_{k}-j)\pi i)}}\right.\no\\
&&\qquad\times\left.\prod_{k=m+1}^{n} \frac{\sh{\frac{\pi}{\rho}
(\beta_{m}-\beta_{k}+(\Lambda_{m}+\Lambda_{k}-j)\pi i)}}{\sh{\frac{\pi}{\rho}
(\beta_{m}-\beta_{k}-(\Lambda_{m}+\Lambda_{k}-j)\pi i)}} \right)
^{\binom{n+l-j-2}{n-1}},
\no
\ena
where $\binom ab$ is the usual binomial coefficient.

We consider the following function
\bea
&&E_{l}(\beta_1,\cdots,\beta_n)=\no\\
&&\left(e^{\mu\sum_{m=1}^{n}\beta_m}\right)^{\binom{n+l-1}n}
\prod_{j=0}^{l-1}\left(\prod_{1\le m<m'\le n}
\frac{S_{2}(i(\beta_m-\beta_{m'})+(\Lambda_{m}+\Lambda_{m'}-j)\pi)}
{S_{2}(i(\beta_m-\beta_{m'})-(\Lambda_{m}+\Lambda_{m'}-j)\pi)}\right)
^{\binom{n+l-j-2}{n-1}}.\label{eq:FUNE}
\ena
By using (\ref{eq:Sprop}), we can check that $E_{l}(\beta_1,\cdots,\beta_n)$
satisfies both (\ref{eq:1de}) and (\ref{eq:2de}). Therefore, we have
\begin{prop}\label{4.1}
\bea
D_{l}(\beta_1,\cdots,\beta_n)
=c_{l}(\rho,\lambda,\mu;\Lambda_1,\cdots,\Lambda_n)
E_{l}(\beta_1, \cdots , \beta_n), \label{eq:prop4.1} 
\ena
where $c_{l}(\rho,\lambda,\mu;\Lambda_1,\cdots,\Lambda_n)$
is a constant independent of $\beta_1,\ldots,\beta_n$.
\end{prop}

In order to determine $c_{l}(\rho,\lambda,\mu;\Lambda_1,\cdots,\Lambda_n)$,
we consider the asymptotics of $D_{l}/E_{l}$ as
\bea
\beta_1, \cdots , \beta_n \in {\bf R}, \quad \beta_1 \ll \cdots \ll \beta_n.
\label{eq:lim}
\ena
Hereafter we use the notation $\sim$ as follows:
\bea
f(\beta_1, \cdots , \beta_n) \sim g(\beta_1, \cdots \beta_n) \quad  \stackrel{{\rm def}}{\Longleftrightarrow} \quad \frac{f(\beta_1, \cdots , \beta_n )}{g(\beta_1, \cdots , \beta_n)} \to 1 \quad {\rm in \, the \, limit }.
\no
\ena

We use the abbreviation $\beta_{mm'}=\beta_m-\beta_{m'}$.
From (\ref{eq:Saasym}), we find
\bea
&&E_{l}(\beta_1,\cdots,\beta_n)\sim
\frac{d_l}{\prod_{L\in{\cal Z}^n_l}P_L(\beta_1,\ldots,\beta_n)}\\
\no
\ena
where
\bea
&&d_l=\exp{\Bigl(\frac{\rho+\lambda}{\rho\lambda}\pi^2i
\left(\binom{n+l-1}n(n-1)\sum_{m=1}^{n}\Lambda_{m}
-\binom n2\binom{n+l-1}{n+1}\right)\Bigr)},\\
&&P_L(\beta_1,\ldots,\beta_n)
=\exp{\Bigl(\frac{2\pi^2}{\rho\lambda}\left(\sum_{m,m'=1}^n
l_m\Lambda_{m'}|\beta_{mm'}|-\sum_{m<m'}l_ml_{m'}\beta_{m'm}
\right)-\mu\sum_{m=1}^{n}l_{m}\beta_{m}\Bigr)}.
\no
\ena
Here we used the equalities
\bea
\sum_{L\in{\cal Z}^n_l}l_{m_1}\cdots l_{m_k}=\binom{n+l-1}{n+k-1}
\quad(1\leq m_1<\cdots<m_k\leq n).
\no
\ena
We have
\bea
\frac{D_{l}(\beta_1,\cdots,\beta_n)}{E_{l}(\beta_1,\cdots,\beta_n)}\sim
\frac1{d_l}\det{\left[P_{L}I(w_{L'}^{(\rho)},w_{L}^{(\lambda)})\right]}
_{L,L'\in{\cal Z}_{l}^{n}}.\label{eq:dfrae} 
\ena

In the following, we consider the asymptotics of
$P_{L}I(w_{L'}^{(\rho)},w_{L}^{(\lambda)})$. We note that  
\bea
I(w_{L'}^{(\rho)},w_{L}^{(\lambda)})=\frac{1}{l!}\sum_{\sigma\in
{\bf S}_{l}}({\rm sgn}\,\sigma)
I(g_{L'}^{(\rho)}(\alpha_{\sigma_1},\ldots,\alpha_{\sigma_l}),
g_{L}^{(\lambda)}(\alpha_1,\ldots,\alpha_l)).
\no
\ena
By changing the integral variables $\alpha_p\to\alpha_p+\beta_{\gamma^{L}(p)}$,
we have
\bea
&&I(g_{L'}^{(\rho)}(\alpha_{\sigma_1},\ldots,\alpha_{\sigma_l}),
g_{L}^{(\lambda)}(\alpha_1,\ldots,\alpha_l))
=\int_C\cdots\int_Cd\alpha_{1}\cdots d\alpha_{l}e^{\mu\sum_{j=1}^{l}
\alpha_{j}+\mu\sum_{m=1}^{n}l_{m}\beta_{m}}\no\\
&&\times\prod_{j=1}^{l}\phi(\alpha_j+\beta_{\gamma^{L}(j)})
\prod_{1\le j<j'\le l}\psi(\alpha_j-\alpha_{j'}
+\beta_{\gamma^{L}(j)\gamma^{L}(j')})\no\\
&&\times g_{L'}^{(\rho)}(\alpha_{\sigma_1}+\beta_{\gamma^{L}(\sigma_1)},\ldots,
\alpha_{\sigma_l}+\beta_{\gamma^{L}(\sigma_l)})
g_{L}^{(\lambda)}(\alpha_1+\beta_{\gamma^{L}(1)},\ldots,
\alpha_l+\beta_{\gamma^{L}(l)}).\label{eq:integrand}
\ena
We set 
\bea
J_{L, L'}^{\sigma}(\alpha_1,\ldots,\alpha_l;\beta_1,\ldots,\beta_n)
=P_{L}(\beta_1,\ldots,\beta_n)
\times({\rm the\,integrand\,of}\,(\ref{eq:integrand})).\no
\ena
From (\ref{eq:Sasym}), we have
\bea
\varphi(x;\Lambda)\sim
\exp{\left(\mp\pi\frac{\rho+\lambda+2\Lambda\pi}{\rho\lambda}x\right)},
\quad(x\rightarrow\infty,\pm{\rm Re}\,x>0)
\label{eq:phiasym1}
\ena
Therefore, if $|{\rm Im}\,x|<K$, we have an estimate
\bea
|\varphi(x;\Lambda)|\leq
\gamma_K\exp{\left(-\pi\frac{\rho+\lambda+2\Lambda\pi}{\rho\lambda}|x|\right)}
\label{eq:phiasym2}
\ena
where $\gamma_K$ is a constant independent of $x$.
Set $\xi(x)=x+|x|$.
Using (\ref{eq:def1}), and (\ref{eq:phiasym2}),
we obtain the following uniform estimate in the asymptotic region
(\ref{eq:lim}).
\bea
&&|J_{L, L'}^{\sigma}(\alpha_1,\ldots,\alpha_l;\beta_1,\ldots,\beta_n)|\le
\gamma\exp\Bigl(\mu\sum_{j=1}^l\alpha_j
+\frac{2\pi^2}{\rho\lambda}\sum_{j<j'}|\alpha_j-\alpha_{j'}|\no\\
&&-\frac{2\pi^2}{\rho\lambda}\sum_{j=1}^l\sum_{m=1}^n\Lambda_m|\alpha_j|
-\frac\pi\lambda\sum_{j=1}^l\xi(\alpha_j)-\frac\pi\rho\sum_{j=1}^l
\xi(\alpha_{\sigma_j}+\beta_{\gamma^L(\sigma_j)\gamma^{L'}(j)})\Bigr).
\label{eq:int}
\ena
Here $\gamma$ is a constant independent of $\alpha_1,\ldots,\alpha_l,
\beta_1,\ldots,\beta_n$.

Now, it is easy to see that $J_{L, L'}^{\sigma}$
is uniformly integrable in (\ref{eq:lim}) if
\bea
\frac{2\pi^2}{\rho\lambda}\left(l-1-\sum_{m=1}^{n}\Lambda_{m}\right)
<\mu<\frac{2\pi}{\lambda}-\frac{2\pi^2}{\rho\lambda}
\left(l-1-\sum_{m=1}^{n}\Lambda_{m}\right).\label{eq:intcon}
\ena
For simplicity, in the following calculation, we assume that
$\rho$ and $\lambda$ are sufficiently large. We have, in particular,
that the region of convergence (\ref{eq:intcon}) is not void.

If $L\not\preceq L'$, then for any $\sigma\in{\bf S}_l$ there exists
$j$ such that $\gamma^L(\sigma_j)>\gamma^{L'}(j)$.
Therefore, by Lebesgue's convergence theorem, we have
\bea
L \not\preceq L' \Longrightarrow P_{L}
I(w_{L'}^{(\rho)}, w_{L}^{(\lambda)}) \to 0 \no
\ena
in the limit (\ref{eq:lim}). Therefore, we get
\bea
\frac{D_{l}(\beta_1,\cdots,\beta_n)}{E_{l}(\beta_1,\cdots,\beta_n)}\sim
\frac1{d_l}
\prod_{L\in{\cal Z}_{l}^{n}}P_{L}I(w_{L}^{(\rho)},w_{L}^{(\lambda)}).
\label{eq:diag}
\ena

Let us calculate the asymptotics of $P_{L}I(w_{L}^{(\rho)},w_{L}^{(\lambda)})$.
Note that
\bea
\gamma^L(\sigma_p)\leq\gamma^L(p)\hbox{ for all $p$}\quad
\Longleftrightarrow\quad\sigma\in{\bf S}_L{\buildrel{\rm def}\over=}
{\bf S}_{l_1}\times\cdots\times{\bf S}_{l_{n}}\subset{\bf S}_{l}.
\no
\ena  
Therefore, we have
\bea
P_{L}I(w_{L}^{(\rho)},w_{L}^{(\lambda)})\sim\frac{1}{l!}
\sum_{\sigma\in{\bf S}_L}
({\rm sgn}\,\sigma)\int\!\cdots\!\int\!d\alpha_1\cdots d\alpha_{l}\,
J_{L,L}^{\sigma}(\alpha_1,\ldots,\alpha_l;\beta_1,\ldots,\beta_n).\no
\ena
Using (\ref{eq:def1}), and (\ref{eq:phiasym1}),
we can calculate the limit of
$J_{L,L}^{\sigma}\,(\sigma\in{\bf S}_{l_1}\times\cdots\times
{\bf S}_{l_{n}})$ as (\ref{eq:lim}):
\bea
P_{L}I(w_{L}^{(\rho)},w_{L}^{(\lambda)})&\to&
\frac{l_1!\cdots l_n!(q^{(\rho)}q^{(\lambda)})^{\sum_{m<m'}l_ml_{m'}}}
{4^{l(n-1)+\sum_{1\leq m<m'\leq n}l_ml_{m'}}l!}
\no\\
&\times&\prod_{m=1}^{n} F_{l_m}^{\Lambda_{m}} \! \! 
\Bigl(\mu +\frac{2 \pi^2}{\rho \lambda} \left( \sum_{j=1}^{m-1}
(l_j -\Lambda_{j})-\sum_{j=m+1}^{n}(l_j-\Lambda_j) \right) \!
-\frac{\rho +\lambda}{\rho \lambda} \pi \Bigr). \label{eq:asymD}
\ena
Here $F_{k}^{\Lambda}(x)$ is given by
\bea
F_{k}^{\Lambda}(x) &=& \int_{C}\cdots\int_{C} d\alpha_1 \cdots d\alpha_{k}
\, e^{x\sum_{j=1}^{k} \alpha_{j}} \prod_{j=1}^{k} \varphi(\alpha_{j};\Lambda)
\! \! \prod_{1 \le j<j' \le k} \! \! \psi(\alpha_{j}-\alpha_{j'})\no\\
&&\qquad{}\times {\rm Skew}\left(  \prod_{1 \le j<j' \le k}
\sh{\frac{\pi}{\rho}(\alpha_{j'}-\alpha_{j}-\pi i)} \right) 
\prod_{1 \le j<j' \le k} \sh{\frac{\pi}{\lambda}
(\alpha_{j'}-\alpha_{j}-\pi i)}. \label{eq:defF}
\ena

From (\ref{eq:prop4.1}), (\ref{eq:diag}), and (\ref{eq:asymD}), we get

\begin{prop}\label{4.2}
\bea 
&&c_{l}(\rho,\lambda,\mu;\Lambda_1,\cdots,\Lambda_n)=
\frac{(q^{(\rho)}q^{(\lambda)})^{\binom n2\binom{n+l-1}{n+1}
+\binom{n+l-1}n\sum_{m=1}^n\Lambda_m}}
{4^{n(n-1)\binom{n+l-1}n+\binom n2\binom{n+l-1}{n+1}}
(l!)^{\binom{n+l-1}{n-1}}}
\prod_{j=1}^{l}(j!)^{n\binom{n+l-j-2}{n-2}}\no\\
&&\quad{}\times\prod_{L\in{\cal Z}_{l}^{n}}
\prod_{m=1}^{n}F_{l_m}^{\Lambda_{m}}\!\!
\Bigl(\mu+\frac{2\pi^2}{\rho\lambda}\left(\sum_{j=1}^{m-1}
(l_j-\Lambda_{j})-\sum_{j=m+1}^{n}(l_j-\Lambda_j)\right)
\!-\frac{\rho+\lambda}{\rho\lambda}\pi\Bigr),\label{eq:prop4.2}
\ena
where $F_{l_m}^{\Lambda_{m}}(x)$ is defined by $(\ref{eq:defF})$.
\end{prop}

\section{Calculation of the integral}

In this section, we find an explicit formula of $F_{l}^{\Lambda}(x)$. 

The integral (\ref{eq:defF}) is absolutely convergent if
\bea
|{\rm Re}\,x|<\frac\pi\rho+\frac\pi\lambda-\frac{2\pi^2}{\rho\lambda}
(l-1-\Lambda).
\label{eq:REGCON}
\ena
In fact, $F^\Lambda_l(x)$ is analytically continued in $x$ to the whole
complex plane. To see this, we derive difference equations satisfied by
$F_{l}^{\Lambda}(x)$. Namely, we relate the values of 
$F_{l}^{\Lambda}$ at $x+\frac{2\pi}\rho$ and $x+\frac{2\pi}\lambda$
to $F_{l}^{\Lambda}(x)$.

\begin{prop}\label{5.1}
The function $F_{l}^{\Lambda}(x)$ satisfies
\bea
\frac{F_{l}^{\Lambda}(x+\frac{\pi}{\lambda})}
{F_{l}^{\Lambda}(x-\frac{\pi}{\lambda})}
=\prod_{k=0}^{l-1}\frac{\ch{\left(\frac{\rho i}{2}x-\frac{\pi^2 i}
{\lambda}(k-\Lambda)\right)}}{\ch{\left(\frac{\rho i}{2}x
+\frac{\pi^2 i}{\lambda}(k-\Lambda)\right)}},\label{eq:equ1}\\
\frac{F_{l}^{\Lambda}(x+\frac{\pi}{\rho})}
{F_{l}^{\Lambda}(x-\frac{\pi}{\rho})}
=\prod_{k=0}^{l-1}\frac{\ch{\left(\frac{\lambda i}{2}x-\frac{\pi^2 i}
{\rho}(k-\Lambda)\right)}}{\ch{\left( \frac{\lambda i}{2}x+\frac{\pi^2 i}
{\rho}(k-\Lambda) \right)}}. \label{eq:equ2}
\ena 
\end{prop}

\begin{proof}
Because of the symmetry between $\rho$ and $\lambda$, it is enough to show
(\ref{eq:equ1}).

We set
\bea
f^{\Lambda}(x | \alpha_1, \cdots , \alpha_l ) &=& e^{x\sum_{j=1}^{l} \alpha_{j}} \prod_{j=1}^{l} \varphi(\alpha_{j};\Lambda) \! \! \! \prod_{1 \le j<j' \le l} \! \! \! \psi(\alpha_{j}-\alpha_{j'})\no\\
& & \quad {}\times {\rm Skew} \! \left(  \prod_{1 \le j<j' \le l} \! \! \sh{\frac{\pi}{\rho}(\alpha_{j'}-\alpha_{j}-\pi i)} \right) \! \!, \\
h(\alpha_1, \cdots , \alpha_l ) &=& \prod_{1 \le j<j' \le l} \! \! \sh{\frac{\pi}{\lambda}(\alpha_{j'}-\alpha_{j}-\pi i)}.
\no
\ena
The poles of $f^\Lambda$ in the variable $\alpha_j$ that are lying
above the contour $C$ are at
$-\Lambda\pi i+\rho i{\bf Z}_{\geq0}+\lambda i{\bf Z}_{\geq0}$
and $\alpha_{j'}+\pi i+\rho i{\bf Z}_{\geq0}+\lambda i{\bf Z}_{\geq0}$.

Consider
\bea
I_k&=&\int_{C}\cdots\int_{C}\left(\int_{C+\rho i}-\int_{C}\right)d\alpha_1
\cdots d\alpha_{l}f^{\Lambda}(x|\alpha_1,\cdots,\alpha_l)\no\\
&&\quad{}\times h(\alpha_1,\cdots,\alpha_l)\sh{\frac{\pi}{\lambda}(\alpha_l
+\Lambda\pi i)}e^{-\frac{\pi}{\lambda}(\sum_{j=1}^{k-1}\alpha_j
-\sum_{j=k}^{l-1}\alpha_j)}.\label{eq:cause}
\ena
The integrand has no poles inside the strip $0<{\rm Im}\,\alpha_l<\rho$,
and the integral is absolutely convergent if $|{\rm Re}\,x|$ is
sufficiently small. Therefore, we have
\bea
I_k&=&0.
\no
\ena

Now we transform the variable $\alpha_l$ to $\alpha_l+\rho i$
so that the contour $C+\rho i$ is modified to $C$.
From (\ref{eq:defphi}) and (\ref{eq:causep3}), we have
\bea
\frac{f^{\Lambda}(x|\alpha_1,\cdots,\alpha_l+\rho i)}
{f^{\Lambda}(x|\alpha_1,\cdots,\alpha_l)}=(-1)^{l}e^{x\rho i}
\frac{\sh{\frac{\pi}{\lambda}(\alpha_l-\Lambda\pi i)}}
{\sh{\frac{\pi}{\lambda}(\alpha_l+\Lambda\pi i+\rho i)}}
\prod_{j=1}^{l-1}\frac{\sh{\frac{\pi}{\lambda}(\alpha_j-\alpha_l-\pi i)}}
{\sh{\frac{\pi}{\lambda}(\alpha_l-\alpha_j-\pi i+\rho i)}}.\label{eq:frule}
\ena
Therefore, the equality $I_k=0$ gives rise to
\bea
0&=&\int_{C}\cdots\int_{C}d\alpha_1\cdots d\alpha_{l}f^{\Lambda}(x|\alpha_1,
\cdots,\alpha_l)e^{-\frac{\pi}{\lambda}(\sum_{j=1}^{k-1}\alpha_j
-\sum_{j=k}^{l-1}\alpha_j)}\no\\
&\times&\left\{(-1)^{l}e^{x\rho i}\sh{\frac{\pi}
{\lambda}(\alpha_l-\Lambda\pi i)}h(\alpha_l,\alpha_1,\cdots,\alpha_{l-1})
-\sh{\frac{\pi}{\lambda}(\alpha_l+\Lambda\pi i)}h(\alpha_1,\cdots,\alpha_l)
\right\}.\label{eq:cause2}
\ena

We use the notation
\bea
a_j=e^{\frac{2 \pi}{\lambda}\alpha_j},
\quad \tau=e^{\frac{\Lambda \pi^2 i}{\lambda}},
\quad q=q^{(\lambda)}.
\no
\ena
Symmetrizing the integrand of (\ref{eq:cause2}), we get
\bea
0=\int_{C}\cdots\int_{C}d\alpha_1 \cdots d\alpha_{l} f^{\Lambda}(x-\frac{\pi}{\lambda} | \alpha_1, \cdots , \alpha_l ) A_{k}(x | \alpha_1, \cdots , \alpha_l), \label{eq:cause3}
\ena
where
\bea
&&A_{k}(x|\alpha_1,\cdots,\alpha_l)\no\\
&&\quad{}={\rm Skew}\left\{\left(-e^{x\rho i}(\tau^{-1}a_1-\tau)a_{k+1}
\cdots a_l-(\tau a_l-\tau^{-1})a_k\cdots a_{l-1}\right)
h(\alpha_1,\cdots,\alpha_l)\right\}.
\no
\ena

Since $(a_j-q^{-2}a_{j+1})h(\alpha_1,\cdots,\alpha_l)$
is symmetric with respect to $\alpha_j$ and $\alpha_{j+1}$, we have
\bea
{\rm Skew}\left( a_j h(\alpha_1, \cdots , \alpha_l ) \right)=q^{-2}{\rm Skew}\left( a_{j+1} h(\alpha_1, \cdots , \alpha_l ) \right). \label{eq:imp}
\ena
Using (\ref{eq:imp}) repeatedly, we get
\bea
{\rm Skew}\left( a_1 a_{k+1} a_{k+2} \cdots a_{l} h(\alpha_1, \cdots , \alpha_l) \right) &=& q^{-2(k-1)}{\rm Skew}\left( a_{k} a_{k+1} \cdots a_{l} h(\alpha_1, \cdots , \alpha_l) \right), \no \\
{\rm Skew}\left( a_{k} \cdots a_{l-1} h(\alpha_1, \cdots , \alpha_l) \right) &=& q^{-2(l-k)} {\rm Skew}\left( a_{k+1} \cdots a_{l} h(\alpha_1, \cdots , \alpha_l) \right),
\no
\ena
and therefore
\bea
A_{k}(x|\alpha_1,\cdots,\alpha_l)&=&-(\tau^{-1}e^{x\rho i}q^{-2(k-1)}+\tau)
{\rm Skew}\left(a_{k}a_{k+1}\cdots a_{l}h(\alpha_1,\cdots,\alpha_l)\right)
\no\\
&&{}+(\tau e^{x\rho i}+\tau^{-1}q^{-2(l-k)}){\rm Skew}\left(a_{k+1}\cdots
a_{l}h(\alpha_1,\cdots,\alpha_l)\right).\label{eq:cause4}
\ena

Substituting (\ref{eq:cause4}) for (\ref{eq:cause3}), we get
\bea
&&(\tau^{-1}e^{x\rho i}q^{-2(k-1)}+\tau)\int_{C}\!\cdots\int_{C}\!\!d\alpha_1
\cdots d\alpha_{l}f^{\Lambda}(x-\frac{\pi}{\lambda}|\alpha){\rm Skew}
\left(a_{k}\cdots a_{l}h(\alpha_1,\cdots,\alpha_l)\right)\no\\
&=&(\tau e^{x\rho i}+\tau^{-1}q^{-2(l-k)})\int_{C}\!\cdots\int_{C}\!\!
d\alpha_1\cdots d\alpha_{l}f^{\Lambda}(x-\frac{\pi}{\lambda}|\alpha)
{\rm Skew}\left(a_{k+1}\cdots a_{l}h(\alpha_1,\cdots,\alpha_l)\right),
\no\\\label{eq:cause5}
\ena
and therefore,
\bea
&&\prod_{k=0}^{l-1}(\tau^{-1}e^{x\rho i}q^{-2k}+\tau)\int_{C}\cdots\int_{C}
\!\!d\alpha_1\cdots d\alpha_{l}f^{\Lambda}(x-\frac{\pi}{\lambda}|\alpha)
{\rm Skew}\left(a_{1}\cdots a_{l}h(\alpha_1,\cdots,\alpha_l)\right)\no\\
&=&\prod_{k=0}^{l-1}(\tau e^{x\rho i}+\tau^{-1}q^{-2k})\int_{C}\cdots\int_{C}
\!\!d\alpha_1\cdots d\alpha_{l}f^{\Lambda}(x-\frac{\pi}{\lambda}|\alpha)
{\rm Skew}\left(h(\alpha_1,\cdots,\alpha_l)\right).
\no
\ena
If $|{\rm Re}\,x|$ is small, both $x\pm\frac\pi\lambda$ lie
in the region (\ref{eq:REGCON}). Therefore, we obtain (\ref{eq:equ1}).
\end{proof}

Consider
\bea
G_{l}^{\Lambda}(x)=\prod_{k=0}^{l-1}\frac{1}{S_{2}\left( \frac{\rho+\lambda}{2} -\pi (k-\Lambda)-\frac{\rho \lambda}{2 \pi}x \right)S_{2}\left( \frac{\rho+\lambda}{2} -\pi (k-\Lambda)+\frac{\rho \lambda}{2 \pi}x \right)} , \label{eq:defG}
\ena
where $S_{2}(x)=S_{2}(x|\rho , \lambda )$. We can check that 
$G_{l}^{\Lambda}$ also satisfies (\ref{eq:equ1}) and (\ref{eq:equ2}).
Therefore, we obtain

\begin{prop}\label{5.2}
\bea
F_{l}^{\Lambda}(x)={\widetilde c}_{l}(\rho, \lambda ; \Lambda )G_{l}^{\Lambda}(x), \label{eq:FCG}
\ena
where ${\widetilde c}_{l}(\rho, \lambda ; \Lambda )$ is
independent of $x$.
\end{prop}

In order to determine ${\widetilde c}_{l}(\rho, \lambda ; \Lambda )$, we consider the asymptotics of $F_{l}^{\Lambda}$ and $G_{l}^{\Lambda}$ as $x \to +i \infty$. From (\ref{eq:Sasym}), we have
\bea
G_{l}^{\Lambda}(x) \sim \exp{\left( -x \pi i \sum_{k=0}^{l-1} (\Lambda -k) \right) }, \quad (x \to +i \infty).
\no
\ena
Therefore,
\bea
{\widetilde c}_{l}(\rho,\lambda;\Lambda)=\lim_{x\to+i\infty}\exp{\left(x\pi i
\sum_{k=0}^{l-1}(\Lambda-k)\right)}F_{l}^{\Lambda}(x).\label{eq:lim0}
\ena

We calculate the limit (\ref{eq:lim0}). 
We use the following equalities.
\bea
&&{\rm Skew}\!\left(\prod_{1\le j<j'\le l}\sh{\frac{\pi}{\rho}(\alpha_{j'}
-\alpha_{j}-\pi i)}\right)=
\frac{[l]_{q^{(\rho)}}!}{l!}
\prod_{1\le j<j'\le l}\!
\!\sh{\frac{\pi}{\rho}(\alpha_{j'}-\alpha_{j})},\label{eq:cause01}\\
&&{e^{\frac{\pi^2i}{2\rho}(l-1)(l-2\Lambda-2)}}
{\rm Skew}\!\left(
e^{-\frac{\pi}{\rho}\sum_{j=1}^{l-1}\alpha_j}
\prod_{j=2}^{l}\sh{\frac{\pi}{\rho}(\alpha_j+\Lambda \pi i)}
\prod_{j=1}^{l-2}\prod_{j'=j+2}^l
\sh{\frac{\pi}{\rho}(\alpha_{j'}-\alpha_{j}-\pi i)} \right)\no\\
&&=\frac1{l!}
\prod_{1\le j<j'\le l}\!\!\sh{\frac{\pi}{\rho}(\alpha_{j'}-\alpha_{j})}.
\label{eq:cause02}
\ena
Combining (\ref{eq:defF}), (\ref{eq:cause01}), and (\ref{eq:cause02}),
we have
\bea
&&\exp{\left(x\pi i\sum_{k=0}^{l-1}(\Lambda-k)\right)}F_{l}^{\Lambda}(x)
=
\frac{[l]_{q^{(\rho)}}![l]_{q^{(\lambda)}}!}{l!}
e^{\frac{\pi^2i}{2\rho}(l-1)(l-2\Lambda-2)}\no\\
&&{}\times\int_C\cdots\int_Cd\alpha_1\cdots
d\alpha_le^{x(\sum_{j=1}^{l}\alpha_j+\pi i\sum_{k=0}^{l-1}
(\Lambda-k))}\prod_{j=1}^{l}\varphi(\alpha_j;\Lambda)\!\!
\prod_{1\le j<j'\le l}\!\!\psi(\alpha_j-\alpha_{j'})\no\\
&&\quad{}\times e^{-\frac{\pi}{\rho}\sum_{j=1}^{l-1}\alpha_j}
\prod_{j=2}^{l}\sh{\frac{\pi}{\rho}(\alpha_j+\Lambda\pi i)}
\prod_{j=1}^{l-2}\prod_{j'=j+2}^l\sh{\frac{\pi}{\rho}
(\alpha_{j'}-\alpha_{j}-\pi i)}
\prod_{1\leq j<j'\leq l}\sh\frac{\pi}{\lambda}
(\alpha_{j'}-\alpha_{j}).\no\\\label{eq:cause03}
\ena

We denote by $I_{l}^{\Lambda}$ the integrand of (\ref{eq:cause03}).
The integral is zero in the limit $x\to+i\infty$ if 
\bea
\Im\left(\sum_{j=1}^{l}\alpha_j+\pi i\sum_{k=0}^{l-1}(\Lambda-k)\right)>0.
\label{eq:cause05}
\ena
Therefore, we can compute the integral by shifting the contour
in the positive imaginary direction, and taking the residues at the poles
which we are crossing.
Note that we assume that $\rho$ and $\lambda$ are large. Therefore,
the only pole we cross is the one at $\alpha_1=\Lambda\pi i$.
Thus we get
\bea
\lim_{x\to i\infty}
\int_C\cdots\int_C d\alpha_1\cdots d\alpha_lI_{l}^{\Lambda}
&=&2\pi i\lim_{x\to i\infty}
\int_{C-\Lambda\pi i+i0}\cdots\int_{C-\Lambda\pi i+i0}d\alpha_2\cdots
d\alpha_l{\rm res}_{\alpha_1=-\Lambda\pi i}I_{l}^{\Lambda}.
\no\\\label{eq:cause04}
\ena
Now the integrand has a pole at $\alpha_2=(1-\Lambda)\pi i$.
Repeating a similar argument, we obtain
\bea
\lim_{x\to i\infty}
\int_C\cdots\int_C d\alpha_1\cdots d\alpha_lI_{l}^{\Lambda}
=(2\pi i)^l{\rm res}_{\alpha_l=(l-1-\Lambda)\pi i}\cdots
{\rm res}_{\alpha_1=-\Lambda\pi i}I_{l}^{\Lambda}.\label{eq:cause06.5}
\ena
From (\ref{eq:lim0}), (\ref{eq:cause03}), (\ref{eq:cause06.5})
(\ref{eq:Sprod}) and (\ref{eq:Sres}), we get

\begin{prop}\label{5.3}
\bea
{\widetilde c}_{l}(\rho,\lambda;\Lambda)
&=&\frac{[l]_{q^{(\rho)}}![l]_{q^{(\lambda)}}!}
{4^{\binom l2}l!}\prod_{k=1}^{l}\frac{S_{2}(\pi)\sqrt{\rho\lambda}}
{S_{2}(k\pi)S_{2}((k-2\Lambda-1)\pi)},\label{eq:ctilde}
\no
\ena
where $S_{2}(x)=S_{2}(x | \rho, \lambda)$. 
\end{prop}

\section{Formula for the determinant}

Combining (\ref{eq:FUNE}), (\ref{eq:prop4.1}), (\ref{eq:prop4.2}),
(\ref{eq:defG}), (\ref{eq:FCG}), and (\ref{eq:ctilde}), we obtain the following
formula for the determinant $D_l(\beta_1,\ldots,\beta_n)$.

\begin{thm}\label{6.1}
\bea
&&\det{\left[I(w_{L'}^{(\rho)},w_{L}^{(\lambda)})\right]
_{L,L'\in{\cal Z}_{l}^{n}}}\no\\
&&{}=\frac
{\left(q^{(\rho)}q^{(\lambda)}\right)^{\binom n2\binom{n+l-1}{n+1}
+\binom{n+l-1}n\sum_{m=1}^n\Lambda_m}
\left(S_{2}(\pi)\sqrt{\rho\lambda}\right)^{n\binom{n+l-1}n}}
{4^{n(n-1)\binom{n+l-1}n+\binom{n+1}2\binom{n+l-1}{n+1}}
(l!)^{\binom{n+l-1}{n-1}}}\no\\&&\times
\prod_{j=1}^{l}\prod_{m=1}^{n}\left(
\frac{[j]_{q^{(\rho)}}!\,[j]_{q^{(\lambda)}}!}
{\prod_{k=1}^{j}S_{2}(k\pi)S_{2}((k-2\Lambda_{m}-1)\pi)}\right)
^{\binom{n+l-j-2}{n-2}}\no\\
&&\quad{}\times\prod_{j=0}^{l-1}\left(\frac
{S_{2}(\frac{\rho\lambda}{2\pi}\mu-(\sum_{m=1}^{n}\Lambda_{m}-j)\pi)}
{S_{2}(\frac{\rho\lambda}{2\pi}\mu+(\sum_{m=1}^{n}\Lambda_{m}-j)\pi)}
\right)^{\binom{n+j-1}{n-1}}E_l(\beta_1,\ldots,\beta_n).
\no
\ena
\end{thm}

\begin{proof}
It remains only to calculate
$c_{l}(\rho,\lambda,\mu;\Lambda_1,\cdots,\Lambda_{m})$
by using (\ref{eq:prop4.2}), (\ref{eq:FUNE}) and (\ref{eq:ctilde}). 

It is easy to see that
\bea
&&\prod_{L\in{\cal Z}_{l}^{n}}
\prod_{m=1}^{n}{\widetilde c}_{l_m}(\rho,\lambda|\Lambda_m)\no\\
&&=\frac
{\left(S_2(\pi)\sqrt{\rho\lambda}\right)^{n\binom{n+l-1}n}}{4^{n\binom{n+l-1}{n+1}}}
\prod_{j=1}^l\prod_{m=1}^n\left(\frac{[j]_{q^{(\rho)}}![j]_{q^{(\lambda)}}!}
{j!\prod_{l=1}^jS_2(l\pi)S_2((l-2\Lambda_m-1)\pi)}\right)
^{\binom{n+l-j-2}{n-2}}.
\no
\ena

From (\ref{eq:defG}), (\ref{eq:Sprop}), and (\ref{eq:Sprod}), we have
\bea
&&G_{l_m}^{\Lambda_m}\Bigl(\mu+\frac{2\pi^2}{\rho\lambda}
\left(\sum_{j=1}^{m-1}(l_j-\Lambda_{j})-\sum_{j=m+1}^{n}(l_j-\Lambda_j)\right)
\!-\frac{\rho+\lambda}{\rho\lambda}\pi\Bigr)\no\\
&&\quad{}=\prod_{l=0}^{l_m-1}\frac
{S_{2}(\frac{\rho\lambda}{2\pi}\mu-\pi(\Lambda^{(m)}-{\widetilde l}_{m}-l))}
{S_{2}(\frac{\rho\lambda}{2\pi}\mu-\pi(\Lambda^{(m-1)}-{\widetilde l}_{m}+l))},
\no
\ena
where $\Lambda^{(m)}=\sum_{j=1}^{m}\Lambda_{j}-\sum_{j=m+1}^{n}\Lambda_{j},$
and $ {\widetilde l}_{m}=\sum_{j=1}^{m-1}l_{j}-\sum_{j=m+1}^{n}l_{j}$.

Now we rewrite 
\bea
\prod_{L\in{\cal Z}_{l}^{n}}\prod_{m=1}^{n}\prod_{k=0}^{l_m-1}\frac
{S_{2}(\frac{\rho\lambda}{2\pi}\mu-(\Lambda^{(m)}-{\widetilde l}_{m}-k)\pi)}
{S_{2}(\frac{\rho\lambda}{2\pi}\mu-(\Lambda^{(m-1)}-{\widetilde l}_{m}+k)\pi)}.
\label{eq:cause000}
\ena

Let us consider the following sets:
\bea
{\cal T}_{\pm}^{(m)}=\bigsqcup_{L\in{\cal Z}_{l}^{n}}
\{{\widetilde l}_{m}\pm k|k=0,1,\cdots,l_m-1\}\quad(m=1,\cdots,n),
\no
\ena
where $\bigsqcup$ means a disjoint union. For $a \in {\bf Z}$, we set
\bea
{\rm mult}_{\pm}^{(m)}(a)=\# \{t \in {\cal T}_{\pm}^{(m)} | t=a \}.
\no
\ena
We set ${\rm mult}_{\pm}^{(0)}(j)={\rm mult}_{\pm}^{(n+1)}(j) \equiv 0$. Then,
we have
\bea
(\ref{eq:cause000})&=&\prod_{m=0}^{n}\prod_{j\in{\bf Z}}
S_{2}\Bigl(\frac{\rho\lambda}{2\pi}\mu-(\Lambda^{(m)}-j)\pi\Bigr)
^{{\rm mult}_{+}^{(m)}(j)-{\rm mult}_{-}^{(m+1)}(j)}\no\\
&=&\prod_{j=0}^{l-1}\left(\frac
{S_{2}\Bigl(\frac{\rho\lambda}{2\pi}\mu-(\Lambda^{(n)}-j)\pi\Bigr)}
{S_{2}\Bigl(\frac{\rho\lambda}{2\pi}\mu-(\Lambda^{(0)}+j)\pi\Bigr)}\right)
^{\binom{n+j-1}{n-1}},\no
\ena
by using
\[
{\rm mult}_{+}^{(n)}(j)={\rm mult}_{-}^{(1)}(-j)=
\cases\binom{n+j-1}{n-1},&\hbox{\rm if }\quad j=0,1,\cdots,l-1;\\
0,&{\rm otherwise},\endcases
\]
\[ {\rm mult}_{+}^{(m)}(j)-{\rm mult}_{-}^{(m+1)}(j)=0 ,  \quad  (m=1, \cdots , n-1), \quad {\rm for \, all} \quad j \in {\bf Z}. \]
This completes the proof.
\end{proof}

{\bf Appendix \quad Double sine function} \\
Here we summarize the property of the double sine function
$S_2(x)=S_{2}(x|\omega_{1},\omega_{2})$ following \cite{JM1}.

We assume $\Re \omega_{1} >0, \Re \omega_{2} >0$. $S_{2}(x|\omega_{1},\omega_{2})$ is a meromorpic function of $x$ and symmetric with respect to $\omega_{1},\omega_{2}$. Its zeros and poles are given by
\bea
{\rm zeros \, at} \> x=\omega_{1}{\bf Z}_{\le 0}+\omega_{2}{\bf Z}_{\le 0}, \quad 
{\rm poles \, at} \> x=\omega_{1}{\bf Z}_{\ge 1}+\omega_{2}{\bf Z}_{\ge 1}. \no
\ena

Its asymptotic behavior is as follows:
\bea
& & \log S_{2}(x)={} \no \\
& & {}=\pm \pi i \! \left( \frac{x^2}{2 \omega_{1} \omega_{2}}-\frac{ \omega_{1}+ \omega_{2}}{2 \omega_{1} \omega_{2}}x-\frac{1}{12}\left(\frac{\omega_{1}}{\omega_{2}}+\frac{\omega_{2}}{\omega_{1}}+3 \right) \right)\! \! + \! o(1), \, (x \to \infty, \, \pm \Im x > 0). \label{eq:Saasym}
\ena
This implies 
\bea
\log S_{2}(a+x)S_{2}(a-x)=\pm \pi i \frac{2a-\omega_{1}-\omega_{2}}{\omega_{1} \omega_{2}}x+o(1), \quad  (x \to \infty, \, \pm \Im x >0). \label{eq:Sasym}
\ena

The following formulae hold:
\bea
\frac{S_{2}(x+\omega_{1})}{S_{2}(x)} &=& \frac{1}{2 \sin \frac{\pi x}{\omega_{2}}}, \label{eq:Sprop} \\
S_{2}(x)S_{2}(-x) &=& -4 {\sin \frac{\pi x}{\omega_{1}}} {\sin \frac{\pi x}{\omega_{2}}},  \label{eq:Sprod} \\
S_{2}(x) &=& \frac{2 \pi}{\sqrt{\omega_{1}\omega_{2}}}x+O(x^2) \quad (x \to 0).
\label{eq:Sres}
\ena

\end{document}